\documentclass[a4paper,11pt]{article}

\usepackage{amsmath,amssymb}

\usepackage{graphicx}
\usepackage{epsfig}
\usepackage{color}

\graphicspath{{Fig/}}

\newenvironment{proof}{\begin{trivlist}
	\item[\noindent]{\it Proof\:}}{\quad $\square$\end{trivlist}}

\newenvironment{exl}{\begin{trivlist}
	\item[\noindent]{\bf Example\:}}{\end{trivlist}}

\newtheorem{dfn}{Definition}[section]
\newtheorem{thm}[dfn]{Theorem}

\newtheorem{lem}[dfn]{Lemma}
\newtheorem{cor}[dfn]{Corollary}

\def\R{{\mathbb R}}

\def\phi{\varphi}
\def\epsilon{\varepsilon}

\def\vol{\mathrm{vol}}

\def\M{{\mathcal M}}
\def\P{{\mathcal P}}

\def\ker{\mathrm{ker}\,}
\def\aff{\mathrm{aff}\,}
\def\span{\mathrm{span}\,}

\title{The Colin de Verdi\`ere number\\ and graphs of polytopes}
\author{Ivan Izmestiev
\thanks{Research for this article was supported by the DFG Research Unit 565 ``Polyhedral Surfaces''.}\\
Institut f\"ur Mathematik\\
Technische Universit\"at Berlin\\
Str. des 17. Juni 136\\
10623 Berlin, Germany\\
{\tt izmestiev@math.tu-berlin.de}}
\date{July 25, 2008}

\begin{document}

\maketitle

\begin{abstract}
The Colin de Verdi\`ere number $\mu(G)$ of a graph $G$ is the maximum corank of a Colin de Verdi\`ere matrix for $G$ (that is, of a Schr\"odinger operator on $G$ with a single negative eigenvalue). In 2001, Lov\'asz gave a construction that associated to every convex 3-polytope a Colin de Verdi\`ere matrix of corank 3 for its 1-skeleton.

We generalize the Lov\'asz construction to higher dimensions by interpreting it as minus the Hessian matrix of the volume of the polar dual. As a corollary, $\mu(G) \ge d$ if $G$ is the 1-skeleton of a convex $d$-polytope.

Determination of the signature of the Hessian of the volume is based on the second Minkowski inequality for mixed volumes and on Bol's condition for equality.
\end{abstract}

\section{Introduction}
\subsection{The Colin de Verdi\`ere number}
At the end of 80's, Yves Colin de Verdi\`ere introduced a graph parameter $\mu(G)$ based on spectral properties of certain matrices associated with the graph~$G$.

\begin{dfn} \label{dfn:Mu}
Let $G$ be a graph with $n$ vertices. A \emph{Colin de Verdi\`ere matrix} for $G$ is a symmetric $n \times n$ matrix $M = (M_{ij})$ with the following properties.
\begin{enumerate}
\item[\rm (M1)] $M$ is a Schr\"odinger operator on $G$, that is
$$
M_{ij} \left\{
\begin{array}{rl}
< 0, & \mbox{ if } ij \mbox{ is an edge of } G;\\
= 0, & \mbox{ if } ij \mbox{ is not an edge of } G \mbox{ and } i \ne j.
\end{array}
\right.
$$
\item[\rm (M2)] $M$ has exactly one negative eigenvalue, and this eigenvalue is simple.
\item[\rm (M3)]
If $X$ is a symmetric $n \times n$ matrix such that $MX = 0$ and $X_{ij} = 0$ whenever $i=j$ or $ij$ is an edge of $G$, then $X=0$.
\end{enumerate}
The set of all Colin de Verdi\`ere matrices for graph $G$ is denoted by $\M_G$. The {\rm Colin de Verdi\`ere number} $\mu(G)$ is defined as the maximum corank of matrices from $\M_G$:
$$
\mu(G) := \max_{M \in \M_G} \dim \ker M.
$$
A Colin de Verdi\`ere matrix of maximum corank is called \emph{optimal}.
\end{dfn}

Basically, the Colin de Verdi\`ere number is the maximum multiplicity of the second least eigenvalue $\lambda_2$ of a discrete Schr\"odinger operator $M$ satisfying a certain stability assumption (M3). By replacing $M$ with $M-\lambda_2 \mathrm{Id}$, we can make the second eigenvalue zero (M2), so that multiplicity becomes corank. Definition \ref{dfn:Mu} was motivated by the study of Schr\"odinger and Laplace operators associated with degenerating families of Riemannian metrics on surfaces.

The parameter $\mu(G)$ turned out to be interesting on its own. In particular, it posesses the minor monotonicity property: if a graph $H$ is a minor of $G$, then $\mu(H) \le \mu(G)$. By the Robertson-Seymour theorem this implies that graphs with $\mu(G) \le n$ can be characterized by a finite set of forbidden minors. For $n$ up to four such characterizations are known and allow nice topological reformulations: e.~g. $\mu(G) \le 3$ iff $G$ is planar (that is doesn't have $K_5$ or $K_{3,3}$ as minors), and $\mu(G) \le 4$ iff $G$ is linklessly embeddable in $\R^3$ (that is doesn't have any graph of the Petersen family as a minor). An overview of results and open problems on the Colin de Verdi\`ere number can be found in \cite{CdV98}, \cite{HLS99}, and \cite{CdV04}. The book \cite{CdV98} deals also with other spectral invariants arising from discrete Schr\"odinger and Laplace operators.

\subsection{Nullspace representations and Steinitz representations}
Let $M$ be a Colin de Verdi\`ere matrix for graph $G$ with $\dim\ker M = d$. Choose a basis $(u_1,\ldots,u_d)$ for $\ker M \subset \R^n$, fix a coordinate system in $\R^n$, and read off the coordinates of $(u_\alpha)$:
$$
(u_1,\ldots,u_d) = (v_1,\ldots, v_n)^\top.
$$
The map that associates to every vertex $i$ of $G$ the vector $v_i \in \R^d$ is called a \emph{nullspace representation} of the graph $G$.

Nullspace representations were studied in \cite{LSc99}. In a subsequent paper \cite{Lov01} Lov\'asz showed that, for a 3-connected planar $G$, the nullspace representation with properly scaled vectors $(v_i)$ realizes $G$ as the skeleton of a convex $3$-polytope. Lov\'asz provided also an inverse construction that associated to every convex 3-polytope with 1-skeleton $G$ a Colin de Verd\`ere matrix of corank~$3$. The proof that the constructed matrix had an appropriate signature was indirect, and a more geometric approach was desirable.

\subsection{Hessian matrix of the volume as a Colin de Verdi\`ere matrix}
In this paper we relate the Lov\'asz construction (that of a matrix from a polytope) to the mixed volumes. Our approach allows a straightforward generalization to higher dimensions. That is, we associate to every $d$-dimensional convex polytope with $1$-skeleton $G$ a Colin de Verdi\`ere matrix for $G$ of corank~$d$.

As a consequence, the graph of a convex $d$-dimensional polytope has Colin de Verdi\`ere number at least $d$. This result is not really new, since it follows from the minor monotonicity of $\mu$, from the fact that the graph of a $d$-polytope has $K_{d+1}$ as a minor \cite{Gru03}, and from $\mu(K_{d+1}) = d$.

Our result is based on the following observation. Take a convex $d$-polytope $P$ and deform it by shifting every facet parallelly to itself. Then the Hessian matrix of the volume of $P$, where partial derivatives are taken with respect to the distances of the shifts, has corank $d$ and exactly one positive eigenvalue. Besides, the mixed partial derivative $\frac{\partial^2\vol(P)}{\partial x_i \partial x_j}$ is positive if the $i$th and the $j$th facets are adjacent, and vanishes otherwise. Thus the negative of the Hessian matrix satisfies conditions (M1) and (M2) from Definition \ref{dfn:Mu}. The condition (M3) follows quite easily, too.

The signature of the Hessian of the volume is encoded in the second Minkowski inequality for mixed volumes together with Bol's characterization of the case of equality. For simple polytopes, the determination of the signature of the Hessian is an essential part in the proof of the Alexandrov-Fenchel inequality.

\subsection{Plan of the paper}
In Section \ref{subsec:LovCon} we recall the Lov\'asz construction of a Colin de Verdi\`ere matrix for the skeleton of a convex 3-polytope $Q$.

After inroducing some terminology and notation in Section \ref{subsec:ParPol}, we show in Section \ref{subsec:HessVol} that the Lov\'asz matrix is minus the Hessian matrix of the volume of the polar dual polytope $Q^*$.

In Section \ref{subsec:InfRig}, dealing with 3-polytopes, we point out an interesting identity (first found and used elsewhere \cite{BI08}) between the Hessian matrix of $\vol(Q^*)$ and the Hessian matrix of another geometric quantity associated with $Q$. This gives another interpretation of the Lov\'asz matrix $M$ and relates the equality $\dim \ker M = 3$ with the infinitesimal rigidity of the polytope $Q$.

In Section \ref{subsec:WhatFails} we discuss the (im)possibility of inverting the construction, that is of finding a convex polytope whose Hessian matrix of the volume equals to a given Colin de Verdi\`ere matrix.

In Section \ref{subsec:NegEig} we give an estimate of the negative eigenvalue (and thus of the spectral gap) for the Hessian matrices of the volume.

Finally, in the Appendix we derive the signature of the Hessian from the second Minkowski inequality and Bol's condition. Although this seems to be a folklore knowledge in narrow circles, we failed to find a written account on this subject.

\subsection{Acknowledgements}
I am grateful to the organizers of the 2006 Oberwolfach conference ``Discrete Differential Geometry'', where the idea of this paper was born. I also thank Ronald Wotzlaw for pointing me out a mistake in a preliminary version.

\section{From a convex polytope to a Colin de Verdi\`ere matrix}
\subsection{Lov\'asz construction}
\label{subsec:LovCon}
Let us recall the Lov\'asz construction of an optimal Colin de Verdi\`ere matrix associated with a polytopal representation of a graph in $\R^3$.

Let $Q \subset \R^3$ be a convex polytope containing the coordinate origin in its interior. Let $G$ be the 1-skeleton of $Q$. We denote the vertices of $G$ by $i,j,\ldots$, and the corresponding vertices of $Q$ by $v_i, v_j, \ldots$. Let $Q^*$ be the polar dual of $Q$. The vertices of $Q^*$ are denoted by $w_f, w_g, \ldots$, where $f,g,\ldots$ are faces of $Q$.

\begin{figure}[ht]
\begin{center}
\input{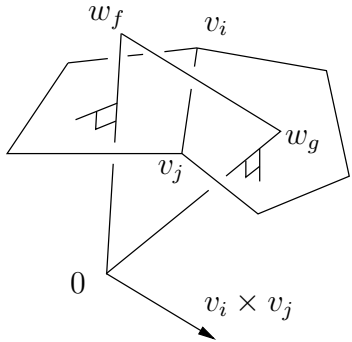}
\end{center}
\caption{To the definition of the matrix $M$.}
\label{fig:PP*}
\end{figure}

For $ij \in G$, consider the edge $v_i v_j$ of $Q$ and the dual edge $w_f w_g$ of $Q^*$, see Figure \ref{fig:PP*}. It is easy to show that the vector $w_f - w_g$ is orthogonal to both vectors $v_i$ and $v_j$, hence parallel to their cross product $v_i \times v_j$. Thus we have
\begin{equation}
\label{eqn:Mij}
w_f - w_g = M_{ij} (v_i \times v_j),
\end{equation}
with $M_{ij} < 0$ (we agree to choose the labeling of $w_f$ and $w_g$ so that we get the correct sign).

Further, consider the vector
$$
v_i' = \sum_j M_{ij} v_j,
$$
where the sum extends over all vertices of $G$ adjacent to $i$. From \eqref{eqn:Mij} it is easy to see that $v_i \times v_i' = 0$. Thus there exists a real number $M_{ii}$ such that
\begin{equation}
\label{eqn:Mii}
v_i' = -M_{ii} v_i.
\end{equation}
Putting $M_{ij} = 0$ for distinct non-adjacent vertices $i$ and $j$ of $G$, we complete the construction of the matrix $M$.

\begin{thm}[Lov\'asz, \cite{Lov01}]
\label{thm:Lov}
The matrix $M$ is a Colin de Verdi\`ere matrix for the graph~$G$.
\end{thm}
The equation \eqref{eqn:Mii} can be rewritten as
\begin{equation}
\label{eqn:Ker}
\sum_j M_{ij} v_j = 0.
\end{equation}
Thus $M$ has corank at least 3. Since $\mu(G) \le 3$ for planar graphs, $M$ is an optimal Colin de Verdi\`ere matrix for $G$.

The proof of Theorem \ref{thm:Lov} goes through a deformation argument, using the fact that the space of convex $3$-polytopes with a given graph is connected.

\subsection{Polytopes with a given set of normals}
\label{subsec:ParPol}
Here we fix some terminology and notation needed in the subsequent sections.

All polytopes in this paper are assumed to be convex. A \emph{facet} of a $d$-dimensional polytope is a $(d-1)$-dimensional face of it.

We will study families of polytopes with fixed facet normals. Let $v_1,\ldots, v_n$ be vectors in $\R^d$ such that the coordinate origin lies in the interior of their convex hull. Consider a $d\times n$ matrix formed by row vectors $v_i^\top$:
$$
V = (v_1,\ldots,v_n)^\top.
$$
\begin{dfn}
Denote by $\P(V)$ the set of all convex polytopes with the outer facet normals $v_1,\ldots, v_n$.
\end{dfn}

Every polytope in $\P(V)$ is the solution set of a system of linear inequalities:
$$
P(x) = \{p \in \R^d\, |\, Vp \le x\},
$$
where $x = (x_i)_{i=1}^n \in \R^n$. Denote by $F_i(x)$ the facet of $P(x)$ with the outer normal $v_i$. We have
$$
F_i(x) = \{p \in P(x)\, |\, v_i^\top p = x_i\}.
$$
The numbers $x_i$ are called the \emph{support parameters} of the polytope $P(x)$. The map $P(x) \mapsto x$ embeds $\P(V)$ into $\R^n$ as an open convex subset. The support parameter $x_i$ is proportional to the signed distance from $0$ to the affine hull of the facet $F_i(x)$:
$$
x_i = \|v_i\| \cdot h_i.
$$

By $\vol_d$ we denote the volume of a $d$-dimensional polytope. We use the subscript because both $\vol_d(P)$ and $\vol_{d-1}(F_i)$ will occur in our formulas. We omit the subscript at $\vol$, when it seems reasonable to do so.

\subsection{Interpreting and generalizing the Lov\'asz construction}
\label{subsec:HessVol}
By definition of the polar dual, we have
$$
Q^* = \{p \in \R^3\, |\, v_i^\top p \le 1 \mbox{ for all }i\}.
$$
Thus $Q^*$ can be viewed as an element of the set $\P(V)$ of polytopes with facet normals $(v_i)_{i \in G}$. In terms of Section \ref{subsec:ParPol}, $Q^* = P(1,\ldots,1)$. Let's vary the support parameters of $Q^*$ and look how does this change its volume.

\begin{lem}
\label{lem:LovVol}
Let $M$ be the matrix constructed in Section \ref{subsec:LovCon}. Then we have
$$
M_{ij} = \left. -\frac{\partial^2 \vol(P(x))}{\partial x_i \partial x_j} \right|_{x = (1,\ldots,1)}
$$
where $P(x)$ is as in Section \ref{subsec:ParPol}.
\end{lem}
\begin{proof}. Let $F_i(x)$ be the facet of $P(x)$ with the normal $v_i$. It is not hard to show that
$$
\frac{\partial \vol_3(P(x))}{\partial x_i} = \frac{\vol_2(F_i(x))}{\|v_i\|}.
$$
Further, for $i \ne j$ we have
$$
\frac{\partial \vol_2(F_i(x))}{\partial x_j} = \frac{\vol_1(F_{ij}(x))}{\|v_j\| \sin \theta_{ij}},
$$
if faces $F_i(x)$ and $F_j(x)$ are adjacent; otherwise this derivative is zero. Here $F_{ij}(x)$ is the common edge of $F_i(x)$ and $F_j(x)$, and $\theta_{ij}$ is the angle between the vectors $v_i$ and $v_j$ (i.~e. the outer dihedral angle at the edge $F_{ij}$). The equations are illustrated in Figure \ref{fig:Der} in one dimension lower and for $\|v_i\| = 1$.

\begin{figure}[ht]
\begin{center}
\input{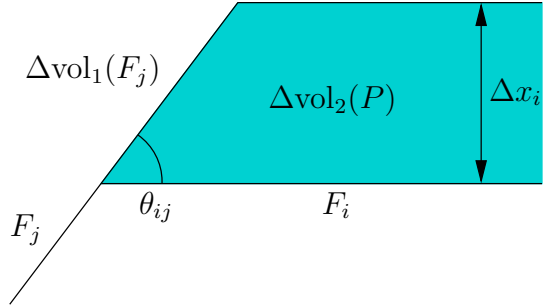}
\end{center}
\caption{Partial derivatives of the volume with respect to the support parameters.}
\label{fig:Der}
\end{figure}

Thus at $x = (1,\ldots,1)$ we have
\begin{equation}
\label{eqn:inej}
\frac{\partial^2 \vol(P(x))}{\partial x_i \partial x_j} = \frac{\vol_1(F_{ij}(x))}{\|v_i\| \|v_j\| \sin \theta_{ij}} = \frac{\|w_f - w_g\|}{\|v_i \times v_j\|} = -M_{ij}
\end{equation}
for all $i \ne j$.

To deal with the case $i=j$, differentiate the well-known identity
$$
\sum_j \vol_2(F_j(x)) \frac{v_j}{\|v_j\|} = 0
$$
with respect to $x_i$. This gives
\begin{equation}
\label{eqn:KerM}
\frac{\partial^2 \vol(P(x))}{\partial x_i^2} v_i + \sum_{j \ne i} \frac{\partial^2 \vol(P(x))}{\partial x_i \partial x_j} v_j = 0.
\end{equation}
In view of \eqref{eqn:Ker} and \eqref{eqn:inej}, we have $\frac{\partial^2 \vol(P(x))}{\partial x_i^2}|_{x=(1,\ldots,1)} = -M_{ii}$.
\end{proof}

Lemma \ref{lem:LovVol} suggests the following generalization of the Lov\'asz construction.

\begin{thm}
\label{thm:HessCdV}
Let
$$
P(x^0) = \{p \in \R^n \, | \, v_i^\top p \le x_i^0 \mbox{ for all }i\}
$$
be a convex polytope with outer facet normals $v_i$ and support parameters $x_i^0, i = 1,\ldots,n$. Let $G$ be the dual 1-skeleton of $P(x^0)$. Then the matrix $M$ defined by
\begin{equation}
\label{eqn:MatrixM}
M_{ij} = \left. -\frac{\partial^2 \vol(P(x))}{\partial x_i \partial x_j} \right|_{x = x^0}
\end{equation}
is a Colin de Verdi\'ere matrix for the graph $G$.

The corank of $M$ is equal to $d$. In particular, $\mu(G) \ge d$ for every graph $G$ that can be realized as the 1-skeleton of a $d$-dimensional polytope.
\end{thm}
\begin{proof}.
Similarly to Lemma \ref{lem:LovVol}, for adjacent facets $F_i$ and $F_j$ we have
$$
\frac{\partial^2 \vol_d(P(x))}{\partial x_i \partial x_j} = \frac{\vol_{d-2}(F_{ij}(x))}{\|v_i\| \|v_j\| \sin \theta_{ij}} > 0,
$$
where $F_{ij}$ is their common $(d-2)$-face, and $\theta_{ij}$ is the angle between $v_i$ and $v_j$. For non-adjacent $F_i$ and $F_j$ this derivative is zero. Therefore matrix $M$ satisfies property (M1) from Definition \ref{dfn:Mu}.

The proof of property (M2) is the most interesting part of the theorem. The signature of the Hessian of the volume is encoded in the second Minkowski inequality for mixed volumes enhanced by Bol's condition for equality.

Theorem \ref{prp:PhiSign} in Section \ref{sec:Append} states in particular that the matrix $M$ has corank $d$. The kernel of $M$ is easy to identify: due to the equation \eqref{eqn:KerM} it consists of the vectors $\xi \in \R^n$ such that $\xi_i = v_i^\top p$ for some vector $p \in \R^d$.

Assuming this description of $\ker M$, let us prove that the matrix $M$ satisfies property (M3). If $MX = 0$, then there are vectors $p_1,\ldots,p_n \in \R^d$ such that $X_{ij} = v_i^\top p_j$ for all $i,j$. Fix $j$. Then by assumption on $X$ we have $p_j \perp v_j$ and $p_j \perp v_i$ for all $ij \in G$. But the normal $v_j$ to the face $F_j$ and the normals to the neighboring faces span the space $\R^d$. Thus we have $p_j = 0$ for all $j$, which implies $X=0$.

As for the last sentence of the theorem, if $G$ is the dual 1-skeleton of a convex polytope $P$, then $G$ is the skeleton of the polar $(P-p)^*$, where $p$ is any interior point of $P$.
\end{proof}

\subsection{Case $d=3$ and infinitesimal rigidity of convex polytopes}
\label{subsec:InfRig}
In the case $d=3$ there is another interpretation of the matrix $M$. As in Section \ref{subsec:LovCon}, let $Q$ be a convex polytope that has skeleton $G$ and contains $0$ in the interior. Triangulate the faces of $Q$ by diagonals and cut $Q$ into pyramids with apices at $0$ and triangles of the triangulation as bases. Denote by $r_i$ the length of the edge that joins $0$ to the vertex $v_i$ of $Q$. Now deform the pyramids by changing the lengths $r_i$ and leaving the lengths of boundary edges constant. During such deformation, the dihedral angles of the pyramids change, and the total angle $\omega_i$ around the $i$-th edge can become different from $2\pi$. By computing the derivatives of $\omega_i$ explicitly, we obtain (\cite{BI08}, Theorem 3.11)
\begin{equation}
\label{eqn:Ide}
\frac{\partial \omega_i}{\partial r_j} = -\frac{\vol_1(F_{ij})}{\sin \theta_{ij}} =  \|v_i\| \|v_j\| \cdot M_{ij},
\end{equation}
where we use the notations from Section \ref{subsec:HessVol}. If we change the variables $x_i$ to $h_i = \|v_i\| \cdot x_i$, so that $h_i$ is the distance of $0$ from $\aff(F_i)$, then the equation \eqref{eqn:Ide} takes a particularly nice form
$$
\frac{\partial \omega_i}{\partial r_j} = -\frac{\partial \vol_2(F_i)}{\partial h_j}.
$$

By \eqref{eqn:Ide}, the matrix $(\frac{\partial \omega_i}{\partial r_j})$ is obtained from the matrix $M$ by multiplying the $i$-th row and the $i$-th column with $\|v_i\|$, for all $i$. This implies
\begin{cor}
The matrix $(\frac{\partial \omega_i}{\partial r_j})$ is an optimal Colin de Verdi\`ere matrix for graph $G$.
\end{cor}

The fact that the matrix $(\frac{\partial \omega_i}{\partial r_j})$ has corank 3 is equivalent to the infinitesimal rigidity of the polytope $Q$. Indeed, every infinitesimal deformation $(dr_i)$ such that $d\omega_i = 0$ for all $i$ gives rise to an infinitesimal isometric deformation of $Q$. The resulting deformation is trivial iff it is produced by moving the apex $0$ inside $Q$.

Another interesting fact is that the matrix $(\frac{\partial \omega_i}{\partial r_j})$ is the Hessian matrix of a geometric quantity related to the polytope $Q$ (deformed by varying $r_i$). Namely, put
$$
S(r) = \sum_{i=1}^n r_i \kappa_i + \sum_{ij \in G} \ell_{ij} \theta_{ij},
$$
where $\kappa_i = 2\pi - \omega_i$ is the ``curvature'' along the $i$-th radial edge, and $\ell_{ij} = \vol_1(F_{ij})$ is the length of the edge $v_i v_j$.
Then the Schl\"afli formula implies
$$
\frac{\partial S}{\partial r_i} = \kappa_i.
$$
Hence
$$
\frac{\partial^2 S(Q)}{\partial r_i \partial r_j} = \frac{\partial^2 \vol(Q^*)}{\partial h_i \partial h_j},
$$
and both matrices are equal to the negative of the Lov\'asz matrix $M$, up to scaling the rows and columns by $\|v_i\|$.

\section{Concluding remarks}

\subsection{What fails in the inverse construction}
\label{subsec:WhatFails}
Let $M$ be a Colin de Verd\`{e}re matrix for the graph $G$. Is there a convex polytope $P$ such that $M$ arises from $P$ as a result of the construction described in Section \ref{subsec:HessVol}? Of course, in general the answer is no, because $G$ must be the dual skeleton of $P$, and $P$ must have dimension $d = \dim\ker M$. In particular, all vertices of $G$ must have degrees at least $d$. But, due to the minor monotonicity of $\mu$, there exist trivalent graphs with $\mu(G)$ arbitrarily large.

Nevertheless, it is worth looking at what fails when we try to reconstruct the polytope $P$ from matrix $M$.

Let $u_1,\ldots, u_d \in \R^n$ be a basis of $\ker M$. Let $v_i^\top$ be the $i$-th row in the matrix $(u_1,\ldots,u_d)$. Then we have
\begin{equation}
\label{eqn:Mink}
\sum_j M_{ij} v_j = 0
\end{equation}
for all $i$. Therefore, the vectors $v_1,\ldots, v_n \in \R^d$ are good candidates for the outer normals to the faces of the polytope $P$. At this point we can already fail, if the following assumptions aren't fulfilled:
\begin{enumerate}
\item
$v_i \ne 0$ for all $i$, and $v_i \ne v_j$ for all $i \ne j$;
\item
for every $i$, the projections $v_{ij}$ of $v_j$ on $v_i^\perp$ for $ij \in G$ satisfy the previous assumption and span $v_i^\perp$.
\end{enumerate}

We proceed assuming that these conditions hold. Codimension $2$ faces $F_{ij}$ of $P$ must be in 1-to-1 correspondence with the edges of $G$, and their volumes are determined by the matrix $M$:
$$
\vol_{d-2}(F_{ij}) = A_{ij} := - M_{ij} \|v_i\| \|v_j\| \sin\theta_{ij},
$$
where $\theta_{ij}$ is the angle between $v_i$ and $v_j$.

\begin{lem}
\label{lem:Facets}
For every $i$, there exists a convex $(d-1)$-dimensional polytope $F_i \subset v_i^\perp$ with outer facet normals $v_{ij}$ and facet volumes $A_{ij}$, $ij \in G$.
\end{lem}
\begin{proof}.
By projecting the equation \eqref{eqn:Mink} on $v_i^\perp$, we obtain
\begin{equation}
\label{eqn:Weight}
\sum_j M_{ij} \cdot v_{ij} = 0.
\end{equation}
Due to $\|v_{ij}\| = \|v_j\| \sin\theta_{ij}$, it follows that
$$
\sum_j A_{ij} \cdot \frac{v_{ij}}{\|v_{ij}\|} = 0.
$$
By Minkowski's theorem \cite[Section 7.1]{Scn93}, this implies the existence of a polytope $F_i$ as stated in the lemma.
\end{proof}

The polytopes $F_i$ in Lemma \ref{lem:Facets} should become facets of the polytope $P$. But here is the second point where the reconstruction can fail: the $j$-th facet $F_{ij}$ of $F_i$ might be different from the $i$-th facet $F_{ji}$ of $F_j$; the only thing we know is $\vol_{d-2}(F_{ij}) = A_{ij} = \vol_{d-2}(F_{ji})$.

In the case $d=3$, however, this suffices: $F_i$ are convex polygons and fit together along their edges to form a polytope $P$. Conditions 1. and 2. above hold if we assume that $G$ is a $3$-connected planar graph \cite{LSc99}. Thus for $3$-connected planar graphs every Colin de Verdi\`ere matrix corresponds to a polytope. This is one of the results of \cite{Lov01}.

The following example shows that even for highly connected graphs the number $\mu(G)$ can be bigger than the maximum dimension of a polytope with 1-skeleton $G$.
\begin{exl}
Let $G_n = K_{2,2,\ldots,2}$ be the multipartite graph on $2n$ vertices (the graph of an $n$-dimensional cross-polytope). By \cite{KLV97}, $\mu(G_n) = 2n-3$ for $n \ge 3$. For $n = 3,4$ the graph $G_n$ can also be represented as the skeleton of a $(2n-3)$-dimensional convex polytope: for $n=3$ this is the octahedron, for $n=4$ the join of two convex quadrilaterals in general position in $\R^5$. For $n \ge 5$, however, there is no $(2n-3)$-dimensional convex polytope with skeleton $G_n$. Indeed, by studying the Gale diagram \cite[Lecture 6]{Zie95} of a $d$-polytope with $d+3$ vertices, one can show that the complement to the graph of such polytope cannot have more than $4$ edges.
\end{exl}

Note that the equation \eqref{eqn:Weight} is reminiscent of the definition of a $(d-2)$-weight in~\cite{McM96}.

\subsection{Negative eigenvalue}
\label{subsec:NegEig}
\begin{thm}
\label{thm:Gap}
Let $\lambda_1$ be the negative eigenvalue of the matrix \eqref{eqn:MatrixM}. Then the following inequality holds:
$$
\lambda_1 \le -d(d-1) \cdot \frac{\vol_d(P(x^0))}{\|x^0\|^2}.
$$
The equality takes place iff
$$
x_i^0 = c \cdot \frac{\vol_{d-1}(F_i(x^0))}{\|v_i\|}
$$
for all $i$ and some constant $c$.
\end{thm}
\begin{proof}.
By induction on $d$, it is easy to show that the function $\vol_d(P(x))$ is a degree $d$ homogeneous polynomial in $x$ as long as the combinatorics of $P(x)$ does not change. For different combinatorics, the polynomials have different coefficients. However, since $\vol_d(P(x))$ is twice differentiable, we can apply Euler's homogeneous function theorem twice at the point $x^0$, independently on how generic the combinatorics of $P(x^0)$ is. This yields
$$
(x^0)^\top M x^0 = -d(d-1) \cdot \vol_d(P(x^0)).
$$

Since $\lambda_1 = \min_{\|\xi\|=1} \xi^\top M \xi$, the inequality follows.

Since $\lambda_1$ is the unique negative eigenvalue of $M$, the inequality turns into equality iff $M x^0 = \lambda x^0$ for some $\lambda$. We have
$$
\sum_j M_{ij} x^0_j = -\frac{1}{\|v_i\|} \sum_j \frac{\partial\vol_{d-1}(F_i(x^0))}{\partial x_j} x_j^0 = -(d-1) \cdot \frac{\vol_{d-1}(F_i(x^0))}{\|v_i\|}.
$$
Thus $M x^0 = \lambda x^0$ is equivalent to $x_i^0 = c \cdot \frac{\vol_{d-1}(F_i(x^0))}{\|v_i\|}$, and the theorem is proved.
\end{proof}

The number $\lambda_2 - \lambda_1$ is called the \emph{spectral gap}. In our case $\lambda_2=0$ by definition. Thus Theorem \ref{thm:Gap} provides an estimate on the spectral gap of the matrix $M$.

Usually, one seeks to make the spectral gap as large as possible, but in order this to make sense for Colin de Verd\`ere matices, one has to choose a matrix norm, \cite[Chapter 5.7]{CdV98}. The norm of the matrix \eqref{eqn:MatrixM} is a function of its coefficients, which have a geometric meaning. Thus, as soon as the choice of a matrix norm is made, one can try to solve the problem of the spectral gap by geometric means (at least for 3-connected planar graphs, for which every optimal Colin de Verdi\`ere matrix can be realized through a polytope).

\begin{appendix}
\section{The second Minkowski inequality for mixed volumes and the signature of the matrix $\left(\frac{\partial^2 \vol}{\partial x_i \partial x_j}\right)$}
\label{sec:Append}

The goal of this appendix is to prove Theorem \ref{prp:PhiSign} that describes the signature of the matrix \eqref{eqn:MatrixM}. The theorem is derived from the second Minkowski inequality for mixed volumes and Bol's condition for equality.

The relation between the theory of mixed volumes and infinitesimal rigidity (as we know, the rank of matrix \eqref{eqn:MatrixM} accounts for the infinitesimal rigidity of the dual polytope, see Section \ref{subsec:InfRig}) was noticed long ago \cite{Bla12, Wey16}. In the decades thereafter this phenomenon seemed to be forgotten. Quite recently, Carl Lee and Paul Filliman \cite{Fil92} discovered it again.

\subsection{The second Minkowski inequality and Bol's condition}

\begin{dfn} \label{dfn:MixVol}
Let $P, Q \subset \R^d$ be convex bodies. A mixed volume of $P$ and $Q$ is a coefficient in the expansion
\begin{equation} \label{eqn:MixVol}
\vol(\lambda P + \mu Q) = \sum_{k=0}^d \binom{d}k \vol(\underbrace{Q,\ldots, Q}_k, \underbrace{P,\ldots, P}_{d-k}) \lambda^{d-k} \mu^k
\end{equation}
with $\lambda, \mu > 0$, where $A+B$ for $A, B \subset \R^d$ denotes the Minkowski sum. In particular,
$$
\vol(P,\ldots, P) = \vol(P).
$$
\end{dfn}
In a similar way one defines the mixed volume of more than two convex bodies. It turns out that the mixed volume is polylinear with respect to the Minkowski addition and multiplication with positive scalars. A proof that the expansion \eqref{eqn:MixVol} takes place and more information on mixed volumes can be found in \cite{Ew96, Scn93}.

\begin{thm}
\label{thm:MinkIneq}
Let $P, Q \subset \R^d$ be convex bodies. Then the following holds:
\begin{enumerate}
\item (The second Minkowski inequality)
\begin{equation} \label{eqn:MinkIneq}
\vol(Q,P,\ldots,P)^2 \ge \vol(P) \cdot \vol(Q,Q,P,\ldots,P).
\end{equation}
\item (Bol's condition)
Assume that $\dim Q = d$. Then equality holds in (\ref{eqn:MinkIneq}) if and only if either $\dim P < d-1$ or $P$ is homothetic to a $(d-2)$-tangential body of $Q$.
\end{enumerate}
\end{thm}

For a proof see \cite[Theorem 6.2.1, Theorem 6.6.18]{Scn93}. Bol's condition was conjectured by Minkowski but proved only decades later by Bol, \cite{Bol43}.

\begin{dfn}
\label{dfn:TanBod}
If $P \subset Q \subset \R^d$ are $d$-dimensional convex polytopes, then $Q$ is called a $p$-tangential body of $P$ iff $P$ has a non-empty intersection with every face of $Q$ of dimension at least $p$.
\end{dfn}

\subsection{Mixed volumes as derivatives of the volume}
By substituting in \eqref{eqn:MixVol} $\lambda = 1$ and $\mu = t$, we obtain
\begin{eqnarray}
\vol(P + tQ) & = & \vol(P) + td\vol(Q,P,\ldots,P)\nonumber\\
&& {}+ t^2\frac{d(d-1)}{2}\vol(Q,Q,P,\ldots,P) + \cdots \label{eqn:VolTay}
\end{eqnarray}
for all $t > 0$, which can be seen as the Taylor expansion of $\vol$. We will look at it in the case when $P$ and $Q$ are polytopes with the same sets of facet normals.

The space $\P(V)$ of all polytopes with outer facet normals $v_1,\ldots,v_n$ is defined in Section \ref{subsec:ParPol}. We want to study the partial derivatives of the volume of $P(x) \in \P(V)$ with respect to the support parameters $x$. For brevity, let's use the notation
$$
\vol(x) := \vol(P(x)).
$$
Similarly, the mixed volume of polytopes from $\P(V)$ will be written as a function of the support parameters:
$$
\vol(x_1,\ldots,x_d) := \vol(P(x_1),\ldots,P(x_d)).
$$

Now we would like to compute $\vol(x+ty)$ with the help of \eqref{eqn:VolTay}. This is not as straightforward as it seems, because the support parameters behave not quite linearly under the Minkowski addition. We have $P(ty) = tP(y)$ for $t>0$. Also we have $P(x) + P(y) \subset P(x+y)$, but the equality doesn't always hold. To describe the cases in which we do have the equality, we need a new definition.

\begin{dfn}
The \emph{normal cone} $N(F,P)$ of the face $F$ of a polytope $P \subset \R^d$ is the set of vectors $w \in \R^d$ such that
$$
\max_{x \in P} (w^\top x) = \max_{x \in F} (w^\top x).
$$
The \emph{normal fan} $N(P)$ is the decomposition of $\R^d$ into the normal cones of the faces of $P$. If the normal fan $N(Q)$ subdivides the normal fan $N(P)$, then we write $N(Q)>N(P)$.
\end{dfn}
Note that the normal fan of a polytope $P \in \P(V)$ has the rays $\R_+v_i$ as 1-dimensional cones. The higher-dimensional cones of the normal fan determine the combinatorics of $P$. Therefore polytopes with equal normal fans are sometimes called strongly isomorphic.

We denote the normal fans of the polytopes from $\P(V)$ by $N(x) := N(P(x))$. The following lemma is classical.

\begin{lem}
\label{lem:NSubd}
If $N(y) > N(x)$, then
$P(x) + P(y) = P(x+y)$.
\end{lem}

Now we are ready to prove
\begin{lem}
\label{lem:VolDer}
Let $y\in \P(V)$ be such that $N(y)>N(x)$. Then
\begin{eqnarray*}
\nabla_y \vol(x) & = & d \cdot \vol(y, x, \ldots, x),\\
\nabla_y^2 \vol(x) & = & d(d-1) \cdot \vol(y, y, x, \ldots, x),
\end{eqnarray*}
where $\nabla_y$ denotes the directional derivative along $y$.
\end{lem}
\begin{proof}.
Due to Lemma \ref{lem:NSubd} we have $P(x+ty) = P(x) + tP(y)$. By substituting $P = P(x)$ and $Q = P(y)$ in \eqref{eqn:VolTay}, we obtain
\begin{eqnarray*}
\vol(x + ty) & = & \vol(x) + td\vol(y,x,\ldots,x)\\
&& {}+ t^2\frac{d(d-1)}{2}\vol(y,y,x,\ldots,x) + \cdots,
\end{eqnarray*}
which implies the lemma.
\end{proof}

\noindent\textbf{Remark.} For polytopes with the same normal fan (``strongly isomorphic polytopes''), there is the following description of mixed volumes. Denote
$$
\P_\Delta(V) = \{P(x) \in \P(V) \,|\, N(x) = \Delta\}.
$$
By induction on $d$, it is easy to show that there exists a homogeneous polynomial $V_\Delta$ of degree $d$ in $n$ variables such that
$$
\vol(P(x)) = V_\Delta(x),
$$
for all $x \in \P_\Delta(U)$. If we use the same symbol $V_\Delta$ to denote the associated symmetric polylinear form, then we have
$$
\vol(P(x^{(1)}),\ldots, P(x^{(d)})) = V_\Delta(x^{(1)},\ldots, x^{(d)})
$$
for all $x^{(1)}, \ldots, x^{(d)} \in \P_\Delta(V)$.

\subsection{From the second Minkowski inequality to the signature of the Hessian of the volume}
By geometric arguments similar to those in the proof of Lemma \ref{lem:LovVol}, the function $\vol$ is twice continuously differentiable on $\P(V)$. Therefore the following definition makes sense.
\begin{dfn} \label{dfn:Phi}
Let $x \in \P(V)$. Define a symmetric bilinear form $\Phi$ on $\R^n$ by
$$
\Phi(\xi,\eta) = \nabla_\eta \nabla_\xi \vol(x).
$$
\end{dfn}

Let $y \in \P(V)$ be such that $N(y) > N(x)$. By combining Euler's homogeneous function theorem and Lemma \ref{lem:VolDer}, we obtain
\begin{eqnarray*}
\Phi(x,x) & = & d(d-1) \, \vol(x, \ldots, x),\\
\Phi(x,y) & = & d(d-1) \, \vol(y, x, \ldots, x),\\
\Phi(y,y) & = & d(d-1) \, \vol(y, y, x, \ldots, x).
\end{eqnarray*}

\begin{lem}
\label{lem:2Sign}
Let $L \subset \R^n$ be a $2$-dimensional vector subspace such that $x \in L$. Then the restriction of the form $\Phi$ to $L$ has signature $(+,-)$ or $(+,0)$.
\end{lem}
\begin{proof}.
Let $y \in \P(V)$ be such that $N(y) > N(x)$. The second Minkowski inequality \eqref{eqn:MinkIneq} applied to $P=P(x)$ and $Q=P(y)$ can be rewritten as
\begin{equation}
\label{eqn:Det}
\det
\begin{pmatrix}
\Phi(x,x) & \Phi(x,y)\\
\Phi(x,y) & \Phi(y,y)
\end{pmatrix}
\le 0.
\end{equation}
Since, moreover, $\Phi(x,x) = d(d-1) \, \vol(P) > 0$, it follows that the restriction of $\Phi$ to $\span\{x,y\}$ has signature $(+,0)$ or $(+,-)$.

It remains to show that every $2$-subspace $L \ni x$ can be represented as $\span\{x,y\}$ with $N(y) > N(x)$. This is true since $x$ is an interior point of the set $\{y \in \P(V) \,|\, N(y) > N(x)\}$. (When we perturb $x$, we can create new faces, but cannot destroy old ones.)
\end{proof}

\begin{lem}
\label{lem:Corank}
The form $\Phi$ has corank $d$.
\end{lem}

\begin{proof}.
Let us exhibit a $d$-dimensional subspace of $\ker \Phi$. Associate with every point $p \in \R^d$ a vector $\overline{p} \in \R^n$ with coordinates
$$
\overline{p}_i = \langle v_i, p \rangle.
$$
The polytope $P(x+\overline{p})$ is the translate of $P(x)$ by $p$. Therefore the directional derivative $\nabla_{\overline{p}} \vol(x)$ vanishes for all $x$, which implies $\Phi(\overline{p}, \eta) = 0$ for all $\eta$. Thus we have
$$
\overline{p} \in \ker\Phi
$$
for all $p \in \R^d$.

Let $\xi \in \ker \Phi$. We need to show that $\xi = \overline{p}$ for some $p \in \R^d$. Denote the span of $x$ and $\xi$ by $L$. Then, by Lemma \ref{lem:2Sign}, the restriction $\Phi|_L$ has signature $(+,0)$ and hence
\begin{equation}
\label{eqn:KerL}
\R\xi = L \cap \ker\Phi.
\end{equation}
Choose $y \in L$ such that $N(y)>N(x)$, and $x$ and $y$ are linearly independent. Then the degeneracy of $\Phi|_L$ means that we have an equality in \eqref{eqn:Det} and thus also in the Minkowski inequality for $P = P(x)$ and $Q = P(y)$. By Bol's condition, see Theorem \ref{thm:MinkIneq}, this happens if and only if the polytope $P(x)$ is homothetic to a $(d-2)$-tangential body of the polytope $P(y)$. By studying Definition \ref{dfn:TanBod}, we see that in $\P(V)$ it is equivalent to $P(x)$ being homothetic to $P(y)$. If $P(x)$ is homothetic to $P(y)$, then $x = \lambda y + \overline{p}$ for some $p \in \R^d$, thus $\overline{p} \in L$. Since $\overline{p} \in \ker\Phi$, it follows that
$$
\R\overline{p} = L \cap \ker\Phi.
$$
By comparing this to \eqref{eqn:KerL}, we conclude that $\xi = \mu \overline{p} = \overline{\mu p}$ for some $\mu \in \R$. Thus the kernel of $\Phi$ is confined to the vectors of the form $\overline{p}$.
\end{proof}

\begin{thm} \label{prp:PhiSign}
The form $\Phi$ has corank $d$ and exactly one positive eigenvalue, which is simple.
\end{thm}

\begin{proof}.
The corank of $\Phi$ is computed in Lemma \ref{lem:Corank}.

The form $\Phi$ has at least one positive eigenvector since $\Phi(x,x) > 0$. Assume that it has more than one. Then there exists a $2$-subspace of $\R^n$ on which $\Phi$ is positively definite. The subgroup of $\mathrm{GL}(\R^n)$ that preserves $\Phi$ acts transitively on the cone of positive directions. Thus there is a positive $2$-subspace $L$ that passes through $x$. This contradicts Lemma \ref{lem:2Sign}. Theorem is proved.
\end{proof}

\end{appendix}

\end{document}